\documentclass[12pt]{article}
\usepackage{tikz}
\usepackage{graphicx}
\usepackage{multirow}
\usepackage{amsmath}
\usepackage{cases}
\usepackage{indentfirst}
\usepackage{mathrsfs}

\input amssym

\newtheorem{lem}{Lemma}[section]%
\newtheorem{thm}[lem]{Theorem}%
\newtheorem{defi}[lem]{Definition}%
\newtheorem{cor}[lem]{Corollary}%
\newtheorem{rem}[lem]{Remark}%
\newenvironment{pf}{\medskip\noindent{\textbf{Proof}.\hspace{0.2cm}}}{\hfill \qed \newline \medskip}

\newcommand{\qed}{\hfill \mbox{\raisebox{0.7ex}{\fbox{}}}}

   \def\d{\delta}

\def\nd{\mathrel{\bigm|\kern-.7em/}}

\begin{document}

\title	{\Large   Non-inner automorphisms of order $p$ in finite $p$-groups admitting cyclic center}
	
	\author{Xuesong Ma$^{1}$, Wei Xu$^{2}$\\ School of Mathematical Sciences, Captial Normal University\\ Beijing, 100048, PRC	
		\thanks{$^1$maxues@cnu.edu.cn; $^2$ xw\_7314@sina.com}}
	
	\maketitle	
	\begin{abstract}
	 Let $G$ be a finite non-abelian $p$-group admitting cyclic center and $p$ be an odd prime. In this paper, we prove that if $C_{G}(Z(\gamma_{3}(G)G^{p}))\nleqslant\gamma_{3}(G)G^{p}$, then $G$ has a non-inner automorphism of order $p$.
	\end{abstract}
	
	{\bf  Keywords: derivation, group cohomology, non-inner automorphism, $p$-group }
	
	{\bf  AMS:20D15, 20D45}	
	
	\section{Introduction}
  It is a longstanding conjecture that every non-abelian $p$-group has a non-inner automorphism of order $p$? The conjecture arose from the Theorem of W. Gasch{\"u}tz, which states that: every non-abelian $p$-group has outer automorphism with order $p$. This conjecture remains an open problem that has not been completely solved.  
  
  P. Schmid proved regular $p$-groups have non-inner automorphism with order $p$ in \cite{peter}; A. Abdollahi proved powerful $p$-groups have non-inner automorphism with order $p$ in \cite{alireza1}; M. T. Benmoussa and Y. Guerboussa proved semi-p-abelian $p$-groups have non-inner automorphism with order $p$ in \cite{mohanmmed}.  
 
   A. Abdollahi proved finite $p$-groups of class $2$ have non-inner automorphisms of order $p$ in \cite{alireza3}. M. Ruscitti, L. Legarreta and M. K. Yadav proved  finite $p$-groups of coclass $3$ have non-inner automorphisms of order $p$ in \cite{rus}.
 
  D. Marian, S. Gheorghe proved strongly Frattinian $p$-groups hold non-inner automorphism with order $p$  in \cite{marian}. A. Abdollahi and S. M. Ghoraishi provided some necessary conditions for a possible counterexample in \cite{mohsen2}.

Throughout this paper, $p$ is an odd prime. Denote $\mathrm{F}_{p}$ as the field of $p$ elements. Let $H$ be a group. By $G$ fixing $H$, we mean that $G$ fixes $H$ elementwise.  $d(G)=\log_{p}(|G/\Phi(G)|)$ is the minimum number of generators of $G$. Let $Z_{i}(G)$ and $\gamma_{i}(G)$ be $i$th term of the upper and lower central series of $G$, respectively.  Denote $\Omega_{i}(G)=\langle g | g^{p^i}=1,g\in G \rangle$ and $G^{p}=\langle g^p | g\in G \rangle$.
 
   Let $G$ be a finite non-abelian $p$-group and $p$ be an odd prime. In this paper, we prove that if $C_{G}(Z(\gamma_{3}(G)G^{p}))\nleqslant\gamma_{3}(G)G^{p}$ and $d(Z(G))=1$, then $G$ has a non-inner automorphism of order $p$. Our main conclusion generalizes the main result of \cite{marian}.

 This paper is organized as follows. To prove our conclusion, in Section 2, we present some basic results about derivations. In Section 3, we present the homological conditions for the existence of derivations that induce isomorphisms of order $p$. And in Section $4$,  we use the existence of some derivations to prove that under certain conditions, non-inner automorphisms of order $p$ exist. In section $5$, we prove Theorem \ref{a2}, which is our main result. 

\section{Preliminary facts and results }

We give some notions and basic results about derivation and automorphism.

Let $H$ be a group and there exists a group action of $G$ on $H$. Denote $C_{H}(G)=\{h\in H|h^{g}=h\;for\;all\;g\in G\}$ as the set of fixed elements of the action of $G$ on $H$.
\begin{defi}\label{2.0}
	A derivation from $G$ to $H$ is defined as a map $\d$ from $G$ to $H$ that satisfies
	\begin{equation*}
		\d(xy)=\d(x)^{y}\d(y)
	\end{equation*}
	for all $x,y\in G$.
\end{defi}	
	
Denote $\mathrm{Der}(G,H)$ as the set of all derivations from $G$ to $H$. And $\mathrm{Der}(G)=\mathrm{Der}(G,G)$, especially. Given $h\in H$, set $\d_{h}(g)=(h^{-1})^{g}h$ for all $g\in G$. Then $\d_{h}$ is a derivation from $G$ to $H$. Let $\mathrm{Ider}(G,H)=\{\d_{h}|h\in H\}$. Say derivations in $\mathrm{Ider}(G,H)$ are inner.  Next, it's a basic fact. 
\begin{lem}
	Let 	\begin{align*}
		\mathcal{M}:\mathrm{Der}(G) &\rightarrow \mathrm{End}(G)\\
		\d&\mapsto \psi
	\end{align*}
where $ \psi(g)=g\d(g)$ for all $g\in G$. Then $\mathcal{M}$ is a bijiective map.
\end{lem}

We say $\psi$ defined the above is endmorphism induced by $\d$. For a normal subgroup $N$ of $G$ and $G$ module $M$, set
\begin{align*}
	\mathrm{Inj}_{N}:\mathrm{Der}(G/N,C_{M}(N))&\rightarrow\mathrm{Der}(G,M)\\
	\d&\mapsto\d_{G}
\end{align*}
that $\d_{G}(g)=\d(gN)$ for all $g\in G$. $\mathrm{Inj}_{N}$ is a injective map. And $\mathrm{Im}\mathrm{Inj}_{N}=\{\d\in\mathrm{Der}(G,M)|\d(g)=0\;for\;all\;g\in N \}$. In the following paper, we use $\mathrm{Der}(G/N,C_{M}(N)) $ instead of $\mathrm{Im}\mathrm{Inj}_{N} $.
To justify whether endmorphism induced by $\d\in \mathrm{Der}(G)$ is an automorphism, we have the following Lemma.
\begin{lem}\label{op}
	Denote $\d\in\mathrm{Der}(G)$. Set
	\begin{align*}
		\xi:G&\rightarrow G/\Phi(G)\\
		g&\mapsto g\d(g)\Phi(G).
	\end{align*}
	If\; $\xi$ is a surjective homomorphism, then endmorphism induced by $\d$ is an automorphism.
\end{lem}
\begin{pf}
	Since $\xi$ is a surjective homomorphism, there exist $g_{1},\cdots,g_{d}$ such that $$\langle g_{1}\d(g_{1}),\cdots,g_{d}\d(g_{d}),\Phi(G)\rangle=G.$$
	Thus $\langle g_{1}\d(g_{1}),\cdots,g_{d}\d(g_{d})\rangle=G$. Endmorphism induced by $\d$ is an  automorphism.
\end{pf}

\begin{rem} When $\d(g)\in\Phi(G)$ for all $g\in G$, map the above is surjective. Then endmorphism induced by $\d$ is an automorphism.
\end{rem}

When we get a $\d\in \mathrm{Der}(G)$, we would justify the order of automorphism inuduced by $\d$. 

	\begin{lem}[Lemma $2.8$ \cite{rus}]\label{porder}
Let $A$ be a normal abelian group of $G$. Let $\d\in \mathrm{Der}(G,A)$ and $\phi$ be an automporphism induced by $\d$. Then \begin{equation*}
			\phi^{n}(g)=\prod\limits_{i=0}\limits^{n}(\d^{i}(g))^{n\choose i}
		\end{equation*}
	for all $g\in G$ where $ \d^{0}(g)=g$ and $\d^{i}(g)=\d(\d^{i-1}(g))$ for all positive integer.
	\end{lem}

Furthermore, if $(\d^{i}(g))^{p}=1$ and $\d^{p}(g)=1$ for all $1\leq i\leq p-1$ and $g\in G$, $\psi$ is automorphism with order $p$. This is the sufficient condition in this paper for us to determine that automorphisms induced by derivations have order $p$.

	Let $N$ be a normal subgroup of $G$ and $\alpha\in \mathrm{Aut}(G)$. Say $\alpha$ fixes $G/N$ if $\alpha(g)\in gN$ for all $g\in G$. And say $\alpha$ fixes $N$ if $\alpha(g)=g$ for all $g\in N$. Denote $\mathrm{Aut}(G/N,N)$ as the set consisting of automorphisms of $G$ that fix both $G/N$ and $N$. Let $ \mathrm{Inn}(G/N,N)$ be a subgroup of $\mathrm{Aut}(G/N,N)$ consisting  of inner automorphisms in $\mathrm{Aut}(G/N,N)$.  In general, there exists a  bijection from $\mathrm{Aut}(G/N,N)$ to $\mathrm{Der}(G/N,Z(N))$. By \cite{peter}, it often happens that $G/N$-mod $Z(N)$ has nontrivial cohomology. Denote $\mathrm{H}^{i}(G/N,Z(N))$ as the $i$th cohomology group of $G/N$ with coefficients in $Z(N)$. When $C_{G}(N)\leq N$, there exists a bijection from $\mathrm{Inn}(G/N,N)$ to $\mathrm{Der}(G/N,Z(N))$, and
	\begin{align*}
	\mathrm{Aut}(G/N,N)/\mathrm{Inn}(G/N,N)&\cong\mathrm{Der}(G/N,Z(N))/\mathrm{Ider}(G/N,Z(N))\\
	&=\mathrm{H}^{1}(G/N,Z(N)).
	\end{align*}

The next two Lemmas aim to show that we only need to consider the module $M$ that satisfies 
$$ \mathrm{H}^{1}(G/N,M)\leq\mathrm{F}_{p}$$

\begin{lem}\label{a1}
	If $C_{G}(\Phi(G))\nleqslant \Phi(G)$, then $G$ has a non-inner automorphism with order $p$.
\end{lem} 
\begin{pf} Prove it by the method of contradiction. According to the main result of \cite{marian}, we only need to consider the case where $C_{G}(Z(\Phi(G)))\leq \Phi(G)$. Then it means that $C_{G}(\Phi(G))\leq C_{G}(Z(\Phi(G)))\leq \Phi(G)$ and it is a contradiction. Then $G$ has a non-inner automorphism with order $p$.
\end{pf}

\begin{lem}\ref{cc}
	 Let $N$ be a normal subgroup of $G$ such that $C_{G}(N)\leq N$. Assume $d(Z(G))<m$. If $\mathrm{H}^{1}(G/N,\Omega_{1}(Z(N)))\cong \mathrm{F}_{p}^{m}$, then $G$ has a non-inner automorphism with order $p$.
\end{lem}
\begin{pf}
	There are $\tau_{1},\cdots,\tau_{m}\in \mathrm{Der}(G/N,\Omega_{1}(Z(N)))\backslash \mathrm{Ider}(G/N,\Omega_{1}(Z(N)))$ such that 
	\begin{equation*}
		\mathrm{H}^{1}(G/N,\Omega_{1}(N)))\cong \oplus_{i=1}^{m}\langle \tau_{i}+\mathrm{Ider}(G/N,\Omega_{1}(Z(N)))\rangle.
	\end{equation*} 

	Let $\eta_{1},\cdots,\eta_{m}$ be the induced automorphisms by $\tau_{1},\cdots,\tau_{m}$. Then $\eta_{1},\cdots,\eta_{m}$ have order $p$.
	
	Suppose $\eta_{1},\cdots,\eta_{m}$ are inner automorphisms. Then there exist $\xi_{1},\cdots,\xi_{m}\in Z(N)\backslash \Omega_{1}(Z(N))$ such that  $\eta_{1},\cdots,\eta_{m}$ are induced by  $\xi_{1},\cdots,\xi_{m}$. Then $\xi_{i}^{p}\in Z(G)$ for all $1\leq i\leq m$.  Let $P=\langle \xi_{1}\rangle\times\cdots\times\langle \xi_{m}\rangle $.
	Since $d(Z(G))< m$,  there exists $a\in Z(G)$ and $t_{1},\cdots,t_{m}$ with $0\leq t_{1},\cdots,t_{m}\leq p-1$ such that $t_{j}\neq 0$ for some $1\leq j\leq m$, and $o(\xi_{1}^{t_{1}}\cdots\xi_{m}^{t_{m}}a)=p$.  So $\xi_{1}^{t_{1}}\cdots\xi_{m}^{t_{m}}a\in \Omega_{1}(Z(N))$. It contradictory to $\xi_{1},\cdots,\xi_{m}\in Z(N)\backslash \Omega_{1}(Z(N))$.

 Then $G$ has a non-inner automorphism with order $p$.
\end{pf}

\begin{lem}\label{kkk}
	Let $N_{1}$ be a normal subgroup of $G$ that satisfies $N\leq N_{1}$, and $M$ be a $\mathrm{F}_{p}(G/N)$-module. Then $ \mathrm{H}^{1}(G/N_{1},C_{M}(N_{1}))\leq \mathrm{H}^{1}(G/N_{1},M).$
\end{lem}

\section{Existence of derivations}
 To simplify description in subsequent Section, we need the following definition.
\begin{defi}
	Let $L$ be a finite $p$-group. $M$ is a $\mathrm{F}_{p}(L)$-module. If $C_{M}(L)\cong\mathrm{F}_{p}$ and $ \mathrm{H}^{1}(L,M)\leq \mathrm{F}_{p}$, we say $M$ is $\mathcal{CR}(L)$-module.
	\end{defi}

 To prove conclusions in Section $4$, we would get some cohomological properties of $\mathcal{CR}(G/\Phi(G))$-module. Specificly, when $ M$ is a $\mathcal{CR}(G/\Phi(G))$-module, we will prove that $$\mathrm{H}^{1}(G/\gamma_{3}(G)G^{p},M)\cong\mathrm{H}^{1}(G/\Phi(G),M)\times \mathrm{F}_{p}^{\mathrm{log}_{p}(|G/\gamma_{3}(G)G^{p}|)}.$$ 

In fact, $\mathrm{H}^{1}(G/\Phi(G),M)$ is trivial.

\begin{lem}\label{3.0}
Assume $L$ is a finite $p$-group. Let $M$ be a $\mathrm{F}_{p}(L)$-module such that  $ C_{M}(L)\cong \mathrm{F}_{p}$ and $\mathrm{H}^{1}(L,M)\cong 1$. Then $M\cong \mathrm{F}_{p}(L)$ as a $ \mathrm{F}_{p}(L)$-module.
	\end{lem}
\begin{pf}
	From Lemma $4$ in \cite{hoe}, $\mathrm{H}^{0}(L,M)\cong 1$. Set map $ \mathrm{T}(m)=m^{\sum\limits_{g\in L}g}$  for all $m\in M$. Since $\mathrm{H}^{0}(L,M)\cong 1$, then $ \mathrm{Im}\mathrm{T}\cong C_{M}(L)$.
	There is $y\in M$ such that $y^{\sum\limits_{g\in L}g}\neq 1$. Set $M_{1}$ is $L$-module generated by $y$. Then $M_{1}\cong \mathrm{F}_{p}(L)$. Notice $\mathrm{H}^{0}(L,\mathrm{F}_{p}(L))\cong 1 $. If $M_{1}$ is a proper submodule of $M$, then there is $z\in M\backslash M_{1}$ such that $ z^{g-1}\in M_{1}$ for all $g\in L$. Since $\mathrm{H}^{0}(L,M_{1})\cong 1$, then $\mathrm{H}^{1}(L,M_{1})\cong 1$. Set $\d$ is a derivation induced by $z$. Thus $\d\in \mathrm{Der}(L,M_{1})\backslash\mathrm{Ider}(L,M_{1})$. It's contradictory to $\mathrm{H}^{1}(L,M_{1})\cong 1$. Hence $ M=M_{1}\cong\mathrm{F}_{p}(L).$
\end{pf}

\begin{lem}\label{pp}
Let $F_{d}$ be the free group generated by the set $X$ which has $d$ elements. Assume that $M$ is a $F_{d}$-module. If $f$ is a map from $X$ to $M$, then there exists $\tau\in \mathrm{Der}(F_{d},M)$ such that $\tau|_{X}=f$.
\end{lem}

\begin{lem}\label{3.3}
	Let $L$ be a finite $p$-group and $L_{1}\leq Z(L)$. $M$ is a $\mathrm{F}_{p}(L)$-module where $L_{1}$ acts on $M$ trivially. If $\mathrm{H}^{1}(L,M)\ncong \mathrm{H}^{1}(L/L_{1},M)$, set $\tau\in \mathrm{Der}(L,M)\backslash\mathrm{Der}(L/L_{1},M)$, then $\tau(y)\in  C_{M}(L)$ for all $y\in L_{1}$.
\end{lem}
\begin{pf}
	Select $\tau\in \mathrm{Der}(L,M)$ and $h\in L_{1}$. Then $\tau([g,h])=0$ for all $g\in L$. From Definition \ref{2.0}, then $\tau([g,h])=\tau(g)^{h-1}+\tau(h)^{1-g}=0$. Thus $\tau(h)\in C_{M}(L)$.
\end{pf}

Throughout this paper, let $D=\langle y_{1},\cdots,y_{d}\rangle$ be a finite $p$-group with $exp(D)=p$, $c(D)=2$, $|D/\Phi(D)|=p^d$ and $|\Phi(D)/\gamma_{3}(D)|=p^{d\choose 2} $. When $d(G)=d$, $G/\gamma_{3}(G)G^{p}$ is isomorphic to a qutoient group of $D$. Then a $G/\gamma_{3}(G)G^{p}$-module is naturally recorded as $ D$-module. Suppose $M$ is a $\mathcal{CR}(G/\Phi(G))$-module, then $M$ is a  
$\mathcal{CR}(D/\Phi(D))$-module.

Denote $J^{1}(\mathrm{F}_{p}(G))=J(\mathrm{F}_{p}(G))$ as the Jacobson radical of group algebra $\mathrm{F}_{p}(G)$. And $J^{i+1}(\mathrm{F}_{p}(G))=J(\mathrm{F}_{p}(G))J^{i}(\mathrm{F}_{p}(G)$ for positive number $i$. Let $K$ be a $\mathrm{F}_{p}(G)$-module. Denote $K_{i}=\{h\in K\;|\;h^{g}=0\;for\;all\;g\in J^{i}(\mathrm{F}_{p}(G)) \}$. Next we construct derivations in $K_{2}$ and $K_{3}$ to prove the main Theorems.

From Definition \ref{2.0}, when we construct derivations, we just give the image of derivations on generators.

In subsequent discussion of the Section, we use additive symbol in computation of abelian group. Suppose $S=\langle a\rangle\cong\mathrm{F}_{p}$. $S$ is regarded as a trivial $D$-module. Then $$\mathrm{H}^{1}(D/\Phi(D),S)\cong \mathrm{Hom}(D,\mathrm{F}_{p})\cong \mathrm{F}_{p}^{d}.$$

There is $\d_{i}\in \mathrm{Der}(D/\Phi(D),S)\backslash\mathrm{Ider}(D/\Phi(D),S)$ such that 
\begin{equation*}
	\d_{i}(y_{j})=\begin{cases}
		a,\quad j=i,\\
		0,\quad j\neq i.
	\end{cases}
\end{equation*}

 Denote $\langle r_{i}\rangle\cong \mathrm{F}_{p}$ for all $1\leq i\leq d$.  For all $g\in D$, set $A=S\oplus_{i=1}^{d}\langle r_{i}\rangle$.
Set
\begin{equation*}
		r_{j}^{g}=r_{j}-\d_{j}(g)
	\end{equation*}. For the other elements, set $(la+\sum\limits_{i=1}\limits^{d}k_{i}r_{i})^{g}=la+\sum\limits_{i=1}\limits^{d}k_{i}r_{i}^{g}$ . Then $A$ is a $\mathrm{F}_{p}(D/\Phi(D))$-module. Since $C_{A}(G)=S$, we have the following result.

\begin{lem}\label{ll} Let $D$ and $A$ be as described above. Let $f$ be a mapping from set $Y=\{y_{1},\cdots,y_{d}\}$ to  $A$. Then there exists $\tau\in \mathrm{Der}(D,A)$ such that $\tau|_{Y}=f$.
\end{lem}
\begin{pf} Assume $X=\{x_1,x_2,\cdots, x_d\}$ and $Y=\{y_1,y_2,\cdots,y_d\}$.
	Let $F_{d}$ be free group defined in Lemma \ref{pp}. Let $ \pi$ be a map from $\mathrm{X}$ to $Y$ such that $\pi(x_{i})=y_{i}$. $\pi$ can induced  a group homomorphism from $F_{d}$ to $D$, and the group homomorphism  is denoted as $\eta$. Then $A$ is a $F_{d}$-module. There exists $\d\in \mathrm{Der}(F_{d},A)$ such that $\d|_{X}=f\pi$. Since $\d(x^{p})=0$ and $\d([x,y,z])=0$ for all $x,y,z\in F_{d}$, we have $\d\in \mathrm{Der}(F_{d}/ker\eta,A)$. And there exists $\tau\in \mathrm{Der}(D,A)$ such that $\tau|_{Y}=f$.
\end{pf}

From Lemma \ref{3.3}, we have the following result.

\begin{lem}\label{00}  Let $D$ and $A$ be as described above.
	 Assume 
	\begin{equation*}
		\begin{aligned}
			\mathrm{Res}:\mathrm{Der}(D,A)&\rightarrow\mathrm{Der}(\Phi(D),S)\\
			\tau&\mapsto \tau|\Phi(D).
		\end{aligned}
		\end{equation*}
	Then $\mathrm{Res}$ is a surjective map.
\end{lem}
\begin{pf}  We construct derivations $\d_{[y_{i},y_{j}]}$, which can be induced by
	\begin{equation*}
		\d_{[y_{i},y_{j}]}(y_{s})=\begin{cases}
			r_{j}, &s=i,\\
			0, & otherwise,
		\end{cases}
	\end{equation*}
	where $1\leq i<j\leq d$. Then
	\begin{equation*}
		\d_{[y_{i},y_{j}]}([y_{s},y_{t}])=\begin{cases}
			-a, &s=i,t=j,\\
			0, &otherwise,
		\end{cases}
	\end{equation*}
	where $1\leq i<j\leq d$, $1\leq s<t\leq d$.

Since $ \mathrm{H}^{1}(\Phi(D),S)=\mathrm{Hom}(\Phi(D),S)\cong\mathrm{F}_{p}^{d\choose2}$, $\mathrm{Res}$ is surjective.
\end{pf}

\begin{lem}\label{ff}
	Let $L$ be a finite $p$-group and $M$ be a $\mathrm{F}_{p}(L)$-module such that  $ C_{M}(L)\cong \mathrm{F}_{p}$. Then there exists a unique submodule $K$ of $\mathrm{F}_{p}(L)$ such that $K\cong M$ as $\mathrm{F}_{p}(L)$-module.
\end{lem}
Since $A\cong \mathrm{F}_{p}(L)_{2}$ as $D$-module, then we have the following result.

\begin{thm}\label{3.2}
	Let $H$ be a maximal subgroup of $D$ and $K$ be a $\mathcal{CR}(D/\Phi(D))$-module. Then $\mathrm{H}^{1}(D,K)\cong\mathrm{H}^{1}(D/\Phi(D),K)\oplus\mathrm{F}_{p}^{d\choose2}$.
\end{thm}
\begin{pf}
	Since $K_{1}\cong \mathrm{F}_{p}$ and $\mathrm{H}^{1}(D/\Phi(D),K)\cong \mathrm{F}_{p}$, we have $K_{2}/K_{1}\cong \mathrm{F}_{p}^{d}$ or $K_{2}/K_{1}\cong \mathrm{F}_{p}^{d-1}$. 
	
	If $K_{2}/K_{1}\cong \mathrm{F}_{p}^{d}$, from Lemma \ref{00}, then $\mathrm{H}^{1}(D,K_{2})\cong \mathrm{H}^{1}(D/\Phi(D),K_{2})\oplus\mathrm{F}_{p}^{d\choose2}$. 
	
	If $K_{2}/K_{1}\cong \mathrm{F}_{p}^{d-1}$, then there are $ \tau_{2},\cdots,\tau_{d}\in K_{2}\backslash K_{1}$ such that $K_{2}/K_{1}=\oplus_{i=2}^{d}\langle \tau_{i}+K_{1} \rangle  $. Assume $K_{1}= \langle a\rangle$. Then there exist minimal generators $h_{1},\cdots,h_{d}$ of $D$ that satisfy:  
	\begin{align*}
		\tau_{i}^{1-h_{j}}=\begin{cases}
			a,\quad &i=j,\\
			0,\quad &i\neq j.
		\end{cases}
	\end{align*}
	From construction in Lemma \ref{00}, $\mathrm{H}^{1}(D,K_{2})\cong \mathrm{H}^{1}(D/\Phi(D),K_{2})\oplus\mathrm{F}_{p}^{d\choose2}$.
	
Let 
	\begin{equation*}
		\begin{aligned}
			\mathrm{Res}_{K}:\mathrm{Der}(D,K)&\rightarrow\mathrm{Der}(\Phi(D),C_{K}(D))\\
			\tau&\mapsto \tau|\Phi(D).
		\end{aligned}
	\end{equation*}
From Lemma \ref{3.3}, $ \mathrm{Res}_{K}$ is well defined.
$ \mathrm{Res}_{K}$ is a surjective map from $\mathrm{H}^{1}(D,K_{2})\cong \mathrm{H}^{1}(D/\Phi(D),K_{2})\oplus\mathrm{F}_{p}^{d\choose2}$. Then $\mathrm{H}^{1}(D,K)\cong\mathrm{H}^{1}(D/\Phi(D),K)\oplus\mathrm{F}_{p}^{d\choose2}$.
\end{pf}

\begin{cor}\label{08}  Assume $d(G)=d$.
Let $K$ be a  $\mathcal{CR}(G/\Phi(G))$-module. Then $$\mathrm{H}^{1}(G/\gamma_{3}(G)G^{p},K)\cong \mathrm{F}_{p}^{\log_{p}(|\Phi(G)/\gamma_{3}(G)G^{p}|)}\oplus\mathrm{H}^{1}(D/\Phi(D),K).$$
\end{cor}
\begin{pf}
	Since $d(G)=d$, there exists $\pi$ such that $\pi$ is a surjective homomorphism from $D$ to $G/\gamma_{3}(G)G^{p}$. Then $K$ is a $\mathcal{CR}(D/\Phi(D))$-module. By Theorem \ref{3.2}, we get $\mathrm{H}^{1}(D/ker\pi,K)=\mathrm{Hom}(\Phi(D)/ker\pi,C_{K}(G))\oplus\mathrm{H}^{1}(D/\Phi(D),K)$. 
	Then $\mathrm{H}^{1}(G/\gamma_{3}(G)G^{p},K)\cong \mathrm{F}_{p}^{\log_{p}(|\Phi(G)/\gamma_{3}(G)G^{p}|)}\oplus\mathrm{H}^{1}(D/\Phi(D),K).$	
	\end{pf}

 \begin{cor}\label{vv}
 	 Assume $d(G)=d$.
 	Let $K$ be a  $\mathcal{CR}(G/\Phi(G))$-module. Let $Q_{1},Q_{2}$ be two subgroups of $G$ such that $\gamma_{3}(G)G^{p}\leq Q_{2}\leq Q_{1}\leq G$ and $ |Q_{1}/Q_{2}|\geq p^{2}$. Then 
 	$$ \mathrm{Der}(G/Q_{2},K)/\mathrm{Der}(G/Q_{1},K)\geq\mathrm{F}_{p}^{2}.$$
 \end{cor}
 \begin{pf}
 	If $ |Q_{1}\Phi(G)/\Phi(G)|\geq p^{2}$, then $$ \mathrm{Hom}(G/\Phi(G),C_{K}(G))/\mathrm{Hom}(G/Q_{1}\Phi(G),C_{K}(G))\geq \mathrm{F}_{p}^{2}.$$ 
 	
 	If $|Q_{1}\Phi(G)/\Phi(G)|= p $, $ \mathrm{Hom}(G/\Phi(G),C_{K}(G))=\mathrm{Hom}(G/Q_{1}\Phi(G),C_{K}(G))\oplus\mathrm{F}_{p}$. And $|Q_{1}\cap\Phi(G)/Q_{2}|\geq p$. From Corllary \ref{08}, $\mathrm{Der}(G/Q_{2},K)\cong\mathrm{Der}(G/(Q_{1}\cap\Phi(G)),K)\oplus\mathrm{F}_{p}^{\mathrm{log}_{p}(|(Q_{1}\cap\Phi(G))/Q_{2}|)} $ .
 	
 	If $Q_{1}\leq \Phi(G)$, from Corllary \ref{08}, $\mathrm{Der}(G/Q_{2},K)\cong\mathrm{Der}(G/Q_{1},K)\oplus\mathrm{F}_{p}^{\mathrm{log}_{p}(|Q_{1}/Q_{2}|)} .$
 	
 	From discussion the above, $\mathrm{Der}(G/Q_{2},K)/\mathrm{Der}(G/Q_{1},K)\geq\mathrm{F}_{p}^{2}.$
 \end{pf}

	Let $M$ be a $\mathrm{F}_{p}(L)$-module and define $K=M\oplus\langle c\rangle$. Choose $\tau\in \mathrm{Der}(L,M)\backslash \mathrm{Ider}(L,M)$. For any $k\in\mathrm{F}_{p},h\in M,g\in L$, we define  $(kc+h)^{g}=c+h^{g}-\tau(g)$. Then, it can be verified that $K$ is also a $\mathrm{F}_{p}(L)$-module. We call $K$ is the module constructed from $\tau$ and $M$.

 We get $\mathrm{H}^{1}(D,K)\cong\mathrm{H}^{1}(D/\Phi(D),K)\oplus\mathrm{F}_{p}^{d\choose2}$ when $ K$ is $\mathcal{CR}(D/\Phi(D))$-module in Theorem \ref{3.2}. Let $1\leq H\leq\Phi(D)$ such that $|\Phi(D)/H|=p^{s}$. Choose $\tau_{1},\cdots,\tau_{s}\in \mathrm{Der}(D,K)\backslash\mathrm{Der}(D/\Phi(D),K)$ such that $$\mathrm{Der}(D/H,K)/\mathrm{Der}(D/\Phi(D),K)=\oplus_{i=1}^{s}\langle \tau_{i}+\mathrm{Der}(D/\Phi(D),K)\rangle .$$ Let $M$ be a $D/H$-module constructed from $\tau_{1},\cdots,\tau_{s}$ and $K$ similar to construction of $A$. 
 \begin{defi}
 	Let $K$ be $\mathcal{CR}(L/\Phi(L))$-module. Take $\gamma_{3}(L)L^{p}\leq L_{1}\leq \Phi(L)$. Say $M$ is $\theta-\mathcal{CR}(L/L_{1})$-module if $M$ has the structure as described above.
 \end{defi}

 Next, we will prove the following Theorem.

\begin{thm}\label{cc}
	Let $K$ be a $\mathcal{CR}(D/\Phi(D))$-module and $M$ be a $\theta-\mathcal{CR}(D)$-module.  Assume 
	\begin{equation*}
		\begin{aligned}
			\mathrm{Res}_{\theta}:\mathrm{Der}(D,M)&\rightarrow\mathrm{Der}(\Phi(D),K_{2})\\
			\tau&\mapsto \tau|\Phi(D).
		\end{aligned}
	\end{equation*}
Then $ \mathrm{Im}\mathrm{Res}_{\theta}\geq \mathrm{F}_{p}^{(d+1).({d\choose2}-1)+1}$ when $ K_{2}\cong \mathrm{F}_{p}(D/\Phi(D))$ as $D/\Phi(D)$-module. And $\mathrm{Im}\mathrm{Res}_{\theta} $ is a surjective map when $ K_{2}\ncong \mathrm{F}_{p}(D/\Phi(D))$ as $D/\Phi(D)$-module.
\end{thm}

In order to prove the above Theorem, we need the following result.
	
	\begin{rem}\label{ee}
Set
\begin{equation*}
	d_{i,j}(y_{s})=\begin{cases}
		r_{j},\quad s=i,\\
		r_{i},\quad s=j,\\
		0,\quad otherwise
	\end{cases}
\end{equation*}
where $1\leq i\leq j\leq d$. Then $d_{i,j}$ can be extended to derivations in $\mathrm{Der}(D/\Phi(D),A) $. Denote $\d_{i,j}$ as derivations extended by $d_{i,j}$.

Suppose $\langle c_{i,j}\rangle\cong \mathrm{F}_{p}$ for all $1\leq i\leq j\leq d$. Let $B=A\oplus_{1\leq i\leq j\leq d}\langle c_{i,j}\rangle$ and
$$ c_{i,j}^{g}=c_{i,j}-\d_{i,j}(g)$$
for all $g\in D$. Then $B$ is a $D$-module.

 And set $\d_{[y_{i},y_{j}]}$ are derivations defined in proof of Lemma \ref{00}. Suppose $\langle c_{[y_{i},y_{j}]}\rangle\cong \mathrm{F}_{p}$ for all $1\leq i<j\leq d$.
Let $R=B\oplus_{1\leq i<j\leq d}\langle c_{[y_{i},y_{j}]}\rangle$ and
$$ c_{[y_{i},y_{j}]}^{g}=c_{[y_{i},y_{j}]}-\d_{[y_{i},y_{j}]}(g)$$
for all $g\in D$.  Then $R$ is a $D$-module and $R\cong\mathrm{F}_{p}(D)_{3}$. And let $d_{j,i}=d_{i,j}$ for $1\leq i\leq j\leq d$. There are derivations of $\mathrm{Der}(D,R)$:
\begin{equation*}
	\d_{[y_{i},y_{j}],t}(y_{s})=\begin{cases}
		d_{t,j},\quad &s=i,\\
		d_{[y_{i},y_{j}]},\quad &s=t,\\
		0,\quad & otherwise.
	\end{cases}
\end{equation*}

\begin{equation*}
	\d_{[y_{i},y_{j}],i}(y_{s})=\begin{cases}
		d_{[y_{i},y_{j}]}+d_{i,j},\quad &s=i,\\
		0,\quad &s\neq i.
	\end{cases}
\end{equation*}

\begin{equation*}
	\d_{[y_{i},y_{j}],j}(y_{s})=\begin{cases}
		d_{j,j},\quad &s=i,\\
		d_{[y_{i},y_{j}]},\quad &s=j,\\
		0,\quad &otherwise.
	\end{cases}
\end{equation*}
And  $\d_{[y_{i},y_{j}],s}([y_{s},y_{t}])=0$, $\d_{[y_{i},y_{j}],s}([y_{i},y_{j}])=r_{s}^{-1}$ where $[y_{s},y_{t}]\neq [y_{i},y_{j}]$. 
Then $\mathrm{H}^{1}(D,R)=\mathrm{H}^{1}(D/\Phi(D),C_{R}(\Phi(D)))\oplus\mathrm{F}_{p}^{d.{d\choose2}}$.
For any maximal subgroup $H$ of $\Phi(D)$  that contains $\gamma_{3}(G)G^p$, $ \mathrm{H}^{1}(D/H,C_{R}(H))=\mathrm{H}^{1}(D/\Phi(D),C_{R}(\Phi(D)))\oplus\mathrm{F}_{p}^{d}$.
\end{rem}

 Next, we have the following Lemma.

\begin{lem}\label{uu}
	Let $ U$ be a $\mathrm{F}_{p}(D)$-module and $W$ be a maximal $\mathrm{F}_{p}(D)$-submodule of $U$. Then 
	$$\mathrm{Der}(D,W)\oplus\mathrm{F}_{p}^{d}\geq\mathrm{Der}(D,U).$$
\end{lem}
\begin{pf}
Choose $ \tau_{1},\cdots,\tau_{d+1}\in\mathrm{Der}(D,U)\backslash\mathrm{Der}(D,W) $. Since $W$ is a maximal $\mathrm{F}_{p}(D)$-submodule of $U$, then $\tau_{i}(h)\in W$ for all $1\leq i\leq d+1$ and $ h\in\Phi(D)$.

 Denote $\mathrm{V}_{i}=\{(k_{1},\cdots,k_{d+1})\in \mathrm{F}_{p}^{d+1}|\sum\limits_{j=1}\limits^{d+1}k_{j}\tau_{j}(y_{i})\in W\}$. Then $\mathrm{dim}\mathrm{V}_{i}\geq d$. Then $ \mathrm{dim}(\mathrm{V}_{1}\cap\mathrm{V}_{2})+d+1\geq \mathrm{dim}(\mathrm{V}_{1})+\mathrm{dim}(\mathrm{V}_{2})$, then $\mathrm{dim}(\mathrm{V}_{1}\cap\mathrm{V}_{2})\geq d-1 $. By induction, we get $ \mathrm{dim}(\cap_{i=1}^{d}\mathrm{V}_{i})\geq1$.
 
 There exists a non-zero vector $(k_{1},\cdots,k_{d+1})\in \mathrm{F}_{p}^{d+1}$ such that $\sum\limits_{i=1}\limits^{d+1}k_{i}\tau_{i}\in \mathrm{Der}(D,W) $.
 Then 
 $\mathrm{Der}(D,W)\oplus\mathrm{F}_{p}^{d}\geq\mathrm{Der}(D,U).$
 \end{pf}

	{\noindent\bf Proof of Theorem \ref{cc}}.
	
\begin{pf}

Take $\tau\in \mathrm{Der}(D,M)\backslash\mathrm{Ider}(D,M)$. Since $\tau([h,g])=0 $ for all $h\in\Phi(D),g\in D$, then $\tau(h)^{g-1}\in K_{1}$ for all $g\in G$. Thus $\tau(h)\in K_{2}$ for all $h\in\Phi(D)$. Hence $\mathrm{Res}_{\theta}$ is well defined.
	
	Since $\mathrm{H}^{1}(D/\Phi(D),K)\leq \mathrm{F}_{p}$, there are three cases: \begin{equation*}
		\begin{aligned}
			&Case 1.\;K\cong \mathrm{F}_{p}(D/\Phi(D))\;as\;D/\Phi(D)-module; \\
			&Case 2.\;K_{2}\;is\;\;isomorphic\;to\;a \;maximal\;D/\Phi(D)-module\;of\; \mathrm{F}_{p}(D/\Phi(D))_{2};\\
			&Case 3.\;K_{2}\;is\;isomorphic\;to\; \;D/\Phi(D)-module\;\mathrm{F}_{p}(D)_{2}.
		\end{aligned}
	\end{equation*}

	Case 1: $M_{3}\cong R$ as $D$ module.
	Then $\mathrm{H}^{1}(D,W)\cong \mathrm{H}^{1}(D/\Phi(D),C_{W}(\Phi(D)))\oplus \mathrm{F}_{p}^{d.{d\choose2}}$. Then  $ \mathrm{Im}\mathrm{Res}_{\theta}\cong \mathrm{F}_{p}^{(d+1).{d\choose2}}$
	
		Case 2: There are $\tau_{2},\cdots,\tau_{d}\in K_{2}\backslash K_{1}$ such that $ K_{2}/K_{1}=\oplus_{i=2}^{d}\langle \overline{\tau_{i}}\rangle$. Set $K_{1}=\langle a \rangle$.
		  There is a minimal generators $g_{1},g_{2},\cdots,g_{d}\in D$ such that
	\begin{equation*} 
		\tau_{i}^{1-g_{j}}=\begin{cases}
			0,\quad i\neq j,\\
			a,\quad i=j.
		\end{cases}
	\end{equation*}
	
	Denote $\tau_{1}\in \mathrm{Der}(D,K_{1})$ that satisfies 
	\begin{equation*} 
		\tau_{1}(g_{j})=\begin{cases}
			0,\quad j\neq 1, \\
			a,\quad j=1.
		\end{cases}
	\end{equation*}
Set $ Q$ is module construct from $\tau_{1} $ and $M$. Set $$ \d_{[g_{1},g_{2}],2},\cdots,\d_{[g_{d-1},g_{d}],2},\cdots,\d_{[g_{d-1},g_{d}],d}$$ are  derivations similar to  defined in Remark\ref{ee}. Then $$ \d_{[g_{1},g_{2}],2},\cdots,\d_{[g_{d-1},g_{d}],2},\cdots,\d_{[g_{d-1},g_{d}],d}\in \mathrm{Der}(D,Q).$$
 Set $$ \d_{1,1},\cdots,\d_{1,d}$$ are  derivations similar to  defined in Remark\ref{ee}. Then $$ \d_{1,1},\cdots,\d_{1,d}\in\mathrm{Der}(D/\Phi(D),C_{Q}(\Phi(D)))\backslash\mathrm{Der}(D,C_{M}(\Phi(D))).$$
Notice $\sum\limits_{i=1}\limits^{d}k_{i}\d_{1,i}\notin \mathrm{Der}(D,C_{M}(\Phi(D)))$ for all non-zero vector $(k_{1},\cdots,k_{d})\in \mathrm{F}_{p}^{d}$. From Lemma\ref{uu}, there exist $ k_{1},\cdots,k_{d}\in\mathrm{F}_{p}^{d}$ such that $\d_{[g_{i},g_{j}],s}+\sum\limits_{i=1}\limits^{d}k_{i}\d_{1,i}\in \mathrm{ Der}(D,M)$. Thus $\mathrm{Im}\mathrm{Res}_{\theta}\cong\mathrm{F}_{p}^{d.{d\choose2}}.$

Cases 3: Set $ \eta\in \mathrm{Der}(D/\Phi(D),K)\backslash\mathrm{Ider}(D/\Phi(D),K)$. Denote $W$ is module constructed from $\eta$ and $M$. Then $W_{3}\cong R$ as $D$ module. Then $\mathrm{H}^{1}(D,W)\cong \mathrm{H}^{1}(D/\Phi(D),C_{W}(\Phi(D)))\oplus \mathrm{F}_{p}^{d.{d\choose2}}$. From Lemma\ref{uu}, $$\mathrm{H}^{1}(D,M)\geq \mathrm{H}^{1}(D/\Phi(D),C_{M}(\Phi(D)))\oplus \mathrm{F}_{p}^{d.({d\choose2}-1)}.$$
	Then  $ \mathrm{Im}\mathrm{Res}_{\theta}\geq \mathrm{F}_{p}^{(d+1).{d\choose2}-d}$.

	From discussion the above, we get conclusion.

\end{pf}
  
 \begin{cor}\label{11}Assume $d(G)=d$.
 	Let $K$ be a $\mathcal{CR}(G/\Phi(G))$ module. Set $\gamma_{3}(G)G^{p}\leq Q_{1}\leq Q_{2}\leq\Phi(G)$ such that $|Q_{2}/Q_{1}|=p$. $M$ is $\theta-\mathcal{CR}(G/Q_{1})$ module. If $\mathrm{log}_{p}(|K_{2}|)=d$, Then 
 	$$ \mathrm{Der}(G/Q_{1},M)/\mathrm{Der}(G/Q_{2},C_{M}(Q_{2}))\geq \mathrm{F}_{p}^{2}.$$
 \end{cor}
\begin{pf}
	Since $d(G)=d$ , there exist surjective homomorphism from $D$ to $G/Q_{1}$. Then $M$ is a  $D$ module. Then there exist $ \theta-\mathcal{CR}(D)$ module $U$ such that $ M$ is $D$ submodule of $U$. Set
	\begin{equation*}
		\begin{aligned}
			\mathrm{Res}_{U}:\mathrm{Der}(D,U)&\rightarrow\mathrm{Der}(\Phi(D),K_{2})\\
			\tau&\mapsto \tau|\Phi(D).
		\end{aligned}
	\end{equation*}
From Theorem\ref{cc}, $ \mathrm{Res}_{U}$ is surjective map. Set $1\leq H\leq\Phi(D)$ that $\pi(H)=Q_{2}/Q_{1} $. Then $\mathrm{Der}(D/ker\pi,M)/\mathrm{Der}(D/H,C_{M}(H))\geq \mathrm{F}_{p}^{2}$. Thus $$ \mathrm{Der}(G/Q_{1},M)/\mathrm{Der}(G/Q_{2},C_{M}(Q_{2}))\geq \mathrm{F}_{p}^{2}.$$

\end{pf}

\begin{cor}\label{oo} Assume $d(G)=d$.
Let $K$ be a $\mathcal{CR}(G/\Phi(G))$ module. Set $\gamma_{3}(G)G^{p}\leq Q_{1}\leq Q_{2}\leq\Phi(G)$ such that $|Q_{2}/Q_{1}|=p^{2}$. $M$ is $\theta-\mathcal{CR}(G/Q_{1})$ module. If $\mathrm{log}_{p}(|K_{2}|)=d+1 $, then there exist a maximal subgroup $H$ of $\Phi(G)$ with $Q_{1}\leq H$ and $ HQ_{2}=\Phi(G)$ such that   $ \mathrm{Der}(G/H,C_{M}(H))/\mathrm{Der}(G/\Phi(G),K)\geq \mathrm{F}_{p}^{2}$
\end{cor}
\begin{pf}
		Since $d(G)=d$ , there exist surjective homomorphism $\pi$ from $D$ to $G/Q_{1}$. Then $M$ is a  $D$ module. Then there exist $ \theta-\mathcal{CR}(D)$ module $U$ such that $ M$ is $D$ submodule of $U$.
		
		Since $\Phi(G)/\gamma_{3}(G)G^{p}$ is elementary abelian $p$-group, there exist two maximal subgroups $E_{1},E_{2} $ of $\Phi(G)$ such that $Q_{1}\leq E_{1}\cap E_{2}$ and $ Q_{2}(E_{1}\cap E_{2})=\Phi(G)$. Denote $\tau_{E_{1}}\in\mathrm{Der}(G/E_{1},K)\backslash\mathrm{Der}(G/\Phi(G),K)$ and $\tau_{E_{2}}\in \mathrm{Der}(G/E_{2},K)\backslash\mathrm{Der}(G/\Phi(G),K) $.
		
		Since $\mathrm{log}_{p}(|K_{2}|)=d+1 $, then $K_{2}\cong\mathrm{F}_{p}(G/\Phi(G))_{2}$ as $G/\Phi(G)$ module. Next, we discuss two cases: 
		$K_{3}\cong\mathrm{F}_{p}(G/\Phi(G))_{3}$ as $G/\Phi(G)$ module and $K_{3}\ncong\mathrm{F}_{p}(G/\Phi(G))_{3}$ as $G/\Phi(G)$ module.
		
		When $K_{3}\cong\mathrm{F}_{p}(G/\Phi(G))_{3}$ as $G/\Phi(G)$ module, there exist non-zero vector $(l_{1},l_{2})\in\mathrm{F}_{p}^{2}$ such that $ l_{1}\tau_{E_{1}}+l_{2}\tau_{E_{2}}\in \mathrm{Ider}(G/Q_{1},M_{3})$. Set $ E_{3}\leq \Phi(G)$ such that $(l_{1}\tau_{E_{1}}+l_{2}\tau_{E_{2}})(x)=0 $ for all $x\in E_{3}$. Notice $E_{3}$ is maximal subgroup of $\Phi(G)$ that $Q_{1}\leq E_{3}$. Denote $E_{4}\leq D$ that $\pi(E_{4})=E_{3}$. Then 
		$C_{M}(E_{3})_{3}\cong C_{R}(E_{4}))_{3}$ as $D/E_{4}$ module. Since $$\mathrm{Der}(D/E_{4},C_{R}(E_{4})_{3})/\mathrm{Der}(D/\Phi(D),C_{R}(\Phi(D)))\geq \mathrm{F}_{p}^{2}  $$
		then $ \mathrm{Der}(G/E_{3},C_{M}(E_{3})_{3})/\mathrm{Der}(G/\Phi(G),K_{3})\geq \mathrm{F}_{p}^{2}$.
	Hence	$$ \mathrm{Der}(G/E_{3},C_{M}(E_{3}))/\mathrm{Der}(G/\Phi(G),K)\geq \mathrm{F}_{p}^{2}.$$
		
		When $K_{3}\ncong\mathrm{F}_{p}(G/\Phi(G))_{3}$ as $G/\Phi(G)$ module, from Remark\ref{ee}, $|K|=p^{\frac{d^{2}+3d}{2}}$. Set $\d\in\mathrm{Der}(G/\Phi(G),K)\backslash\mathrm{Ider}(G/\Phi(G),K) $. Denote $_{1}M$ is $G/Q_{1}$ module that construct from $\d$ and $M$. Since $ |\mathrm{F}_{p}(G/\Phi(G))|=p^{p^{d}}\geq p^{\frac{d^{2}+3d}{2}+2}$, then $$\mathrm{H}^{1}(G/\Phi(G),C_{_{1}(M)}(\Phi(G)))\neq 0.$$ 
		From proof of Lemma\ref{uu}, there exists $\eta\in\mathrm{Der}(G/E_{1},C_{M}(E_{1})_{3})\backslash\mathrm{Ider}(G/E_{1},C_{M}(E_{1})_{3})$ such that $ \eta(h)\neq0$ for $h\in \Phi(G)\backslash E_{1}$. Then $ \mathrm{Der}(G/E_{1},C_{M}(E_{1}))/\mathrm{Der}(G/\Phi(G),K)\geq \mathrm{F}_{p}^{2}$.
		
		We get conclusion.
\end{pf}

	\section{ Automorphisms with order $p$}
In this section, $G$ is finite $p$-group with cyclic center.

	\begin{lem}\label{12}
	 If $\Omega_{1}(Z(G))\nleqslant \gamma_{3}(G)G^{p}$, then $G$ has a non-inner automorphism of order $p$.
	\end{lem}
\begin{pf}For $\Omega_{1}(Z(G))\nleqslant \gamma_{3}(G)G^{p}$, there is $a\in \Omega_{1}(Z(G))\backslash\gamma_{3}(G)G^{p}$. Then there is $b\in Z_{2}(G)\backslash\Phi(G)$ such that $b^{p}\in Z(G)$. Notice that $b\in C_{G}(Z(\Phi(G)))$. By main result in \cite{marian}, $G$ has non-inner automorphism with order $p$.
\end{pf}

If $|\Phi(G)/\gamma_{3}(G)G^{p}|=1$, then $G$ is a powerful $p$-group. $G$ has a non-inner automorphism of order $p$ by \cite{alireza1}.
In the remain part of this section, thus we always assume that $\Omega_{1}(Z(G))\leq \gamma_{3}(G)G^{p}$ and $|\Phi(G)/\gamma_{3}(G)G^{p}|\geq p$.

	 Let $\{g_{1},\cdots,g_{d}\}$ be a generator set of $G$ and
	\begin{equation*}
		G/\Phi(G)=\langle \overline{g_{1}}\rangle\times\cdots\times\langle \overline{g_{d}}\rangle.
	\end{equation*}

Let $|\Phi(G)/\gamma_{3}(G)G^{p}|=p^{T}$. And $P_{1},\cdots,P_{T}$ are subgroups of $\Phi(G)$ that satisfy 
$ \gamma_{3}(G)G^{p}= P_{T}\leq \cdots\leq P_{1}\leq P_{0} =\Phi(G)$ and $|P_{i}/P_{i+1}|=p$ for all $0\leq i\leq T-1$.

	\begin{thm}\label{01}			
Assume $ C_{G}(P_{i})\leq P_{i}$ and all automorphisms induced by derivations in $ Der(G/P_{i},\Omega_{1}(Z(P_{i})))$ are inner.  Then there exists an automorphism $\phi$ with order $p$ such that $\phi|_{P_{i+1}}=\mathrm{id}_{P_{i+1}}$ and $\phi|_{P_{i}}\neq\mathrm{id}_{P_{i}}$.
	\end{thm}

	\begin{pf}	From Lemma \ref{cc}, then there is a $G/P_{i}$-module $W_{1}$ of $\Omega_{1}(Z(P_{i}))$  such that  $\mathrm{H}^{1}(G/P_{i},W_{1})\lesssim \mathrm{F}_{p}$ and $C_{G}(W_{1})\cong\mathrm{F}_{p}$. From Lemma \ref{kkk}, $$\mathrm{H}^{1}(G/\Phi(G),C_{W_{1}}(\Phi(G)))\lesssim \mathrm{F}_{p}.$$ From Theorem\ref{3.2}, there is $\tau\in \mathrm{Der}(G/P_{i+1},W_{1})$ such that $\tau(P_{i})\neq 1$ and $\tau(P_{i+1})=1$. 
		By Lemma \ref{porder}, since $W_{i}\leq \Phi(G)$ and $\Omega_{1}(Z(G))\leq\gamma_{3}(G)G^{p}$,  we have $\tau_{p}(g)=1$ for all $g\in G$. The automorphism induced by $\tau$ is order $p$.
		
	\end{pf}

	\begin{cor}\label{02}
	Assume $ C_{G}(P_{i})\leq P_{i}$ and all automorphisms induced by derivations in $ Der(G/P_{i},\Omega_{1}(Z(P_{i})))$ are inner. If $C_{G}(P_{i+1})\nleqslant P_{i+1}$ and $Z(C_{G}(P_{i+1}))>Z(P_{i+1})$, then $G$ has a non-inner automorphism of order $p$.
	\end{cor}
	\begin{pf} Prove the conclusion by contradiction.
		Suppose all automorphisms of $G$ of order $p$ are inner. For $x\in C_{G}(P_{i+1})$, then inner automorphism induced by $x$ fixes $Z(C_{G}(P_{i+1}))$. By Theorem \ref{01}, there exists an automorphism $\phi$ such that $o(\phi)=p$ and $\phi(g)\neq g$ for arbitrary $g\in G\backslash P_{i+1}$. And $\phi|_{P_{i+1}}=\mathrm{id}_{P_{i+1}}$. It's a contradiction.	\end{pf}

\begin{rem}\label{17}
	From Corllary \ref{02}, if $|C_{G}(P_{i})P_{i}/P_{i}|=p$, then $G$ has a non-inner automorphism of order $p$ .
\end{rem}

\begin{thm}\label{011}
Assume  $ C_{G}(P_{i})\leq P_{i}$ and all autmorphisms induced by derivations in $\mathrm{Der}(G/P_{i},\Omega_{1}(Z(P_{i})))$ are inner. If $ C_{G}(P_{i+1})\nleqslant P_{i+1}$, then $G$ has a non-inner automorphism of order $p$.
\end{thm}	

\begin{pf}	
	Denote $W=\Omega_{1}(Z(P_{i}))$ and $I=\{g\in Z(\Phi(G))|[g,h]\in W\}$. Then $\mathrm{H}^{1}(G/P_{i},W)\leq\mathrm{F}_{p}$.
	
When $C_{G}(P_{i+1})\leq P_{i}$, then $|C_{G}(P_{i+1})P_{i+1}/P_{i+1}|=p$. From Remark \ref{17}, $G$ has a non-inner automorphism with order $p$.

So we suppose $C_{G}(P_{i+1})\nleqslant P_{i}$.   When $ Z(C_{G}(P_{i+1}))>Z(P_{i+1})$, $G$ has a non-inner automorphism with order $p$. 

Therefore, we only need to consider the case: $Z(C_{G}(P_{i+1}))=Z(P_{i+1})$.

We claim that $|C_{G}(P_{i+1})P_{i}/P_{i}|=p$. 

If $|C_{G}(P_{i+1})P_{i}/P_{i}|\geq p^{2}$, since $ C_{G}(P_{i})\leq P_{i}$ and all autmorphisms induced by derivations in $\mathrm{Der}(G/P_{i},\Omega_{1}(Z(P_{i})))$ are inner, then $\mathrm{H}^{1}(G/P_{i},W)\leq\mathrm{F}_{p}$. From Lemma\ref{kkk}, $\mathrm{H}^{1}(G/\Phi(G),C_{W}(\Phi(G)))\leq\mathrm{F}_{p}.$  From Corllary\ref{vv}, $$\mathrm{Der}(G/P_{i},C_{W}(\Phi(G)))/\mathrm{Der}(G/(C_{G}(P_{i+1})P_{i}),C_{W}(\Phi(G)))\geq \mathrm{F}_{p}^{2}.$$  Then $ |Z(P_{i})/C_{Z(P_{i})}(C_{G}(P_{i+1}))|\geq p^{2}$. Notice $Z(P_{i})\cap P_{i+1}\leq Z(P_{i+1})$. Then $Z(P_{i})\cap P_{i+1}\leq C_{Z(P_{i})}(C_{G}(P_{i+1}))$.  Since $|P_{i}/P_{i+1}|=p$, we have $|Z(P_{i})/(Z(P_{i})\cap P_{i+1})|\leq p$. It's contradictory to $$ |Z(P_{i})/C_{Z(P_{i})}(C_{G}(P_{i+1}))|\geq p^{2}.$$

Suppose $|C_{G}(P_{i+1})P_{i}/P_{i}|=p$. When $|C_{G}(P_{i+1})P_{i+1}/P_{i+1}|=p$, from Remark \ref{17}, $G$ has a non-inner automorphism with order $p$. We assume $|C_{G}(P_{i+1})P_{i+1}/P_{i+1}|=p^{2}$. Set $h\in P_{i}\backslash P_{i+1},y\in C_{G}(P_{i+1})\backslash P_{i}$. There exists $$\tau\in \mathrm{Der}(G/P_{i+1},C_{W}(\Phi(G)))\backslash\mathrm{Der}(G/P_{i},C_{W}(\Phi(G)))$$ such that $\tau(h)\in \Omega_{1}(Z(G))$. The automorphism inuduced by $\tau$ is order $p$. Suppose the automorphism is inner. Set $g\in G\backslash P_{i}$ that $[x,g]=\tau(x)$ for all $x\in G $. Notice $g\in C_{G}(P_{i+1})\backslash P_{i}$.

 Since $\mathrm{H}^{1}(G/\Phi(G),C_{W}(\Phi(G)))\leq\mathrm{F}_{p}$, then either $C_{W}(\Phi(G))_{2}\cong\mathrm{F}_{p}^{d} $ or $C_{W}(\Phi(G))_{2}\cong\mathrm{F}_{p}^{d+1} $.

Next, we discuss the case $C_{W}(\Phi(G))_{2}\cong\mathrm{F}_{p}^{d} $. Assume $y\in G\backslash\Phi(G)$. If there exist $ y_{1}\in\Omega_{1}(Z_{2}(G)) $ such that $[y,y_{1}]\in\Omega_{1}(Z(G))\backslash\{1\}$, then there exist $y_{2}\in C_{W}(\Phi(G))$ such that $[y,y_{2}]\in Z_{2}(G)\backslash Z(G)$. Since $y\in C_{G}(P_{i+1})$, then $y_{1},y_{2}\in P_{i}\backslash P_{i+1}$. It's contradictory to $ |P_{i}/P_{i+1}|=p$. Then $[y,x]=1$ for all $ x\in C_{W}(\Phi(G))$. There exist $g_{1}\in Z_{2}(G),a\in Z(G)$ such that $(gg_{1}a)^{p}=1$.
Notice $g=y^{l}z$ that $l\in\mathrm{F}_{p}$ and $z\in P_{i}$. Then $[gg_{1}a,x]=1$ for all $x\in C_{W}(\Phi(G))$.
Set $O=C_{W}(\Phi(G))\times\langle gg_{1}a \rangle$ and $H\leq \Phi(G)$ that $[gg_{1}a,x]=1$ for all $x\in H$. Notice $H$ is a maximal subgroup of $\Phi(G)$ that $\gamma_{3}(G)G^{p}\leq H$. Then there exist $\d\in\mathrm{Der}(G/H,O)\backslash Ider(G/H,O)$ such that $\d(h)\neq1 $. If $\d(gg_{1}a)\in O\backslash W_{1}$, then $[\d(gg_{1}a),h]\neq 1$. Since $[gg_{1}a,W]=1,[h,W]=1$, $\d([gg_{1}a,h])=[\d(gg_{1}a),h][gg_{1}a,\d(h)]\neq 1$. It's contradictory to $\d(H)=1$. Then $\d(gg_{1}a)\in W_{1}$. Thus $\d_{3}(x)=1$ for all $x\in G$.  Automorphism induced by $\d$ is order $p$. Notice the automorphism fix $P_{i+1}$ not $h$. Then $G$ has non-inner automorphism with order $p$.

Assume $y\in \Phi(G)$. Set $M=\{x\in W|[x,z]\in Z_{2}(G)\;for\;all\;z\in \Phi(G)\backslash P_{i}\}$. From Corllary\ref{11}, exist $\d\in \mathrm{Der}(G/P_{i},M)\backslash\mathrm{Ider}(G/P_{i},M)$ such that $\d(y)\neq1$. Then $C_{W}(P_{i})/C_{W}(C_{G}(P_{i+1})P_{i})\geq \mathrm{F}_{p}^{2}$. It's contradictory to $|P_{i}/P_{i+1}|=p$.

Next, we disscuss the case $ C_{W}(\Phi(G))_{2}\cong\mathrm{F}_{p}^{d+1}$.
If $y\notin\Phi(G) $, set $H$ is maximal subgroup of $G$ that $y\notin H$, then $ C_{I}(H).Z(G)/Z(G)\geq\mathrm{F}_{p}^{2}$. It's contradictory to $|P_{i}/P_{i+1}|=p$. We assume $C_{G}(P_{i})\leq \Phi(G)$.

 If $ C_{W}(\Phi(G))\cong \mathrm{F}_{p}(G/\Phi(G))$ as $G/\Phi(G)$ module, when $W\nleqslant P_{i+1}$, $$|C_{W}(P_{i})/C_{W}(C_{G}(P_{i+1})P_{i})|\geq p^{2}.$$It's contradictory to $|P_{i}/P_{i+1}|=p$. Assume $W\leq P_{i+1}$.
  Then exist $0\leq l\leq p-1$ such that $o(yh^{l})=p$.  Set $O=C_{W}(\Phi(G))\times\langle yh^{l} \rangle $. There exists $\eta\in \mathrm{Der}(G/P_{i+1}, O)\mathrm{Ider}(G/P_{i+1}, O) $ such that $\eta(h)\in Z_{2}(G)\backslash Z(G) $. Similarly, the automorphism inuduced by $\eta$ is order $p$. Then $G$ has non-inner automorphism with order $p$.  

If $ C_{W}(\Phi(G))\ncong \mathrm{F}_{p}(G/\Phi(G))$ as $G/\Phi(G)$ module, for $z\in C_{G}(P_{i+1})P_{i+1}\backslash P_{i+1}$, set $H_{z}$ is maximal subgroup of $\Phi(G)$ that $[z,x]=1$ for all $x\in H_{z}$. Notice $Q_{1}\leq H_{z}$.  There exists $c\in C_{G}(P_{i+1})P_{i+1}\backslash P_{i+1}$ such that $ O(c)=p$ and $\mathrm{H}^{1}(G/H_{c}, C_{W}(\Phi(G))\times\langle c \rangle )\geq\mathrm{F}_{p}^{2}$. Set $\beta\in\mathrm{Der}(G/H_{c}, C_{W}(\Phi(G))\times\langle c \rangle )\backslash \mathrm{Ider}(G/H_{c}, C_{W}(\Phi(G))\times\langle c \rangle ) $ that $\beta(x)\neq 1$ for $x\in\Phi(G)\backslash H_{z}$.
 Similarly, automorphism induced by $\beta$ is order $p$.
 $G$ has non-inner automorphism with order $p$.
 
 Now, we get conclusion.

\end{pf}

\section{Main result and its proof}

Now, we prove our main conclusions.

\begin{thm}\label{04}
	Let $G$ be a finite $p$-group with cyclic center where $p$ is an odd prime. If  $ C_{G}(\gamma_{3}(G)G^{p})\nleqslant \gamma_{3}(G)G^{p}$, then $G$ has a non-inner automorphism of order $p$.
\end{thm}	

\begin{pf}	
	From Lemma \ref{a1}, we assume $ C_{G}(\Phi(G))\leq \Phi(G)$. Since $P_{0}=\Phi(G)$, by Theorem \ref{011} , $G$ has a non-inner automorphism with order $p$ when $ C_{G}(\gamma_{3}(G)G^{p})\nleqslant \gamma_{3}(G)G^{p}$. 
\end{pf}

	\begin{thm}\label{a2}
	Let $G$ be a finite $p$-group with cyclic center where $p$ is an odd prime. If $C_{G}(Z(\gamma_{3}(G)G^{p}))\nleqslant\gamma_{3}(G)G^{p}$, then $G$ has a non-inner automorphism of order $p$.
\end{thm}
	
	\begin{pf}
		By Theorem \ref{04}, we need to consider that $C_{G}(\gamma_{3}(G)G^{p})=Z(\gamma_{3}(G)G^{p})$.
		
		Suppose  $h\in C_{G}(Z(\gamma_{3}(G)G^{p}))\backslash \gamma_{3}(G)G^{p}$. From Theorem \ref{01}, $G$ has an automorphism of order $p$ which fixes $\gamma_{3}(G)G^{p}$ and $G/\gamma_{3}(G)G^{p}$  but not $h$. The automorphism is non-inner.
	\end{pf}

\end{document}